\documentclass[11pt]{amsart}
\usepackage{latexsym,amsmath,amsthm,amssymb,enumerate,psfrag,epsfig,graphicx}
\usepackage[latin1]{inputenc}
\usepackage{dsfont,textcomp}
\newcommand{\PT}[1]{\mathbf{P}_{\!#1}} %semi-groupe P_1
\newcommand{\GI}{\mathbf{L}} %generator
\newcommand{\var}{\mathbf{Var}} %variance
\newtheorem{theorem}{Theorem}
\newtheorem{proposition}[theorem]{Proposition}
\newtheorem{lemma}[theorem]{Lemma}
\newtheorem{corollary}[theorem]{Corollary}
\newtheorem{claim}[theorem]{Claim}

\theoremstyle{definition} 

\newtheorem*{assumption}{Hypothesis}

\theoremstyle{remark} 
\newtheorem{remark}[theorem]{Remark}

\begin{document}   %----------------------------------------------------
%%%---------------------------------------------------------------------
\date{\today}
\title{Slow decay of Gibbs measures with heavy tails} 

\author[C. Roberto]{\textbf{\; Cyril Roberto}}

\address{{\bf {Cyril} ROBERTO}\\
Laboratoire d'Analyse et Math\'ematiques Appliqu\'ees- UMR 8050,
Universit\'es de Paris-est, Marne la Vall\'ee \\
5 Boulevard Descartes, Cit\'e Descartes, Champs sur Marne\\
77454 Marne la Vall\'ee Cedex 2, FRANCE.}

\email{cyril.roberto@univ-mlv.fr}

\begin{abstract}
We consider Glauber dynamics reversible with respect to Gibbs measures with heavy tails. 
Spins are unbounded. The interactions are bounded and finite range. The self potential enters into two classes of 
measures, $\kappa$-concave probability measure and sub-exponential laws, for which it is known that
no exponential decay can occur.
We prove, using coercive inequalities, that the associated infinite volume semi-group 
decay to equilibrium polynomially and  stretched exponentially, respectively.
Thus improving and extending previous results by Bobkov and Zegarlinski.
\end{abstract}

\maketitle

%%%%%%%%%%%%%%%%%%%%%%%%%%%%%%%%%%%%%%%%%%%%%%%%%%%%%%%%%%%%%%%%%%%%%%%%%%%%%%%%%%%%%%%%%%%%%%%%%%%%%%%%%%%%%
%%%%%%%%%%%%%%%%%%%%%%%%%%%%%%%%%%%%%%%%%%%%%%%%%%%%%%%%%%%%%%%%%%%%%%%%%%%%%%%%%%%%%%%%%%%%%%%%%%%%%%%%%%%%%
%%%%%%%%%%%%%%%%%%%%%%%%%%%%%%%%%%%%%%%%%%%%%%%%%%%%%%%%%%%%%%%%%%%%%%%%%%%%%%%%%%%%%%%%%%%%%%%%%%%%%%%%%%%%%

\section{Introduction}

In the past decades, the study of functional inequalities deserved a lot of attention, not only on the side
of theoretical probability and analysis, but also in statistical mechanics.
This is due to the numerous fields of application: differential geometry, analysis of p.d.e., 
concentration of measure phenomenon, isoperimetry, trends to equilibrium in deterministic and stochastic evolutions... 

The most popular functional inequalities are Poincar\'e and log Sobolev. Both are now well understood in
many situations. We refer to 
\cite{bakry}, \cite{Dav}, \cite{grossbook}, \cite{guionnet-zegarlinski}, \cite{ledoux}, \cite{martinelli}, \cite{Ro99}, \cite{Wbook},  \cite{ane} for an introduction.

Very recently, generalisations of Poincar\'e and log Sobolev inequalities were introduced and studied
by probabilists and analysts.
Let us mention weak Poincar\'e or super Poincar\'e inequalities, Orlicz-Poincar\'e or Orlicz-Sobolev inequalities, $F$-Sobolev inequalities, weighted Poincar\'e or weighted log Sobolev inequalities, 
modified log-Sobolev inequalities etc. To give a complete picture of the literature is out of reach. See \cite{bobkov-zeg2,DGGW,Wang08,BLweight,cggr,barthe-k,BCR1,barthe-roberto} (and references therein) 
for some of the most recent publications.

Few of those recent advances  have been used so far in statistical mechanics, at the notable exception of 
\cite{bobkov-zeg2,zitt}.

On the other hand, in the statistical mechanics community, progress have been done in the study of 
Poincar\'e and log Sobolev inequalities for large classes of models coming from the physics literature.
Again, to give an updated list of publications is out of reach. Let us mention 
\cite{landim-yau, landim-yau2,bcc,cancrini-cm,boudou,cmrt}. 

This paper intends to use advances in both communities in order to improve and extend some
results of Bobkov and Zegarlinski \cite{bobkov-zeg2} on the decay to equilibrium of some
unbounded spins systems. We believe that the techniques
coming more specifically from one community
should be more largely exploited by the other one. This paper is one step in this direction. 

\bigskip

If a lot of results are known for log-concave probability measures, not so much has been proved
for measures with heavy tails (let us mention 
\cite{rockner-wang,bcr2,BLweight,bob07,bobkov-zeg2,WangOrlicz}). 
In this paper the focus is on such measures with heavy tails (informally with tails larger than exponential) 
and our aim is to prove decay to equilibrium
of unbounded spin systems in infinite dimension.

\bigskip

Now we introduce and discuss one of the main tool we shall use, namely, the weak Poincar\'e inequality.
Consider for simplicity a one dimensional probability measure $d\nu=Z_V^{-1}e^{-V}$ with ($Z_V=\int e^{-V(u)} du<\infty$).
Then, a weak Poincar\'e inequality asserts that 
\begin{equation} \label{eq:intro}
\var_\nu (f) %= \int \left( f -\int f d\nu \right)^2 d\nu 
\leq \beta(s) \int (f')^2 d\nu + s \mathrm{Osc}(f)^2
\qquad \forall f : \mathbb{R} \to \mathbb{R}, \forall s>0
\end{equation}
where $\beta : (0,\infty) \to \mathbb{R}$ is a rate function associated to the weak Poincar\'e inequality and 
$\mathrm{Osc}(f)= \sup f - \inf f$. In case when $\lim_0 \beta < \infty$, then the weak Poincar\'e inequalities
reduce to the usual Poincar\'e inequalities. Most of the information is encoded in the behaviour of 
$\beta$ near the origin. Moreover, note that $\var_\nu(f) \leq \frac{1}{4} \mathrm{Osc}(f)$, 
hence, only the values $s \in (0,1/4)$ are relevant.

Weak Poincar\'e inequalities, in the form  
\eqref{eq:intro}, have been introduced by R\"ockner and Wang \cite{rockner-wang}. However, 
inequalities with a free parapeter
have a long story in analysis, see {\it e.g.} \cite{nash,Dav,liggett,bertini-z1,bertini-z2}.

Using capacity techniques (Hardy type inequalities \cite{mazya,muckenhoupt,BR03})
the best possible rate function for $V_p(x)=|x|^p$, $p \in (0,1)$ was computed in \cite{bcr2}: 
$\beta(s)= c_p \log \left( \frac{2}{s \wedge 2} \right)^{2(\frac{1}{p}-1)}$. Also, in this case, it is known 
(see \cite[Chapter 5]{ane}) that $\nu$ does not satisfies the usual Poincar\'e inequality.
Equivalently there is no exponential decay to equilibrium.
However, by standard differentiation (see Section \ref{sec:wp}), one gets that the semi-group $(\mathbf{S}_t)_{t \geq 0}$
associated to the one dimensional generator $L=\frac{d^2}{du^u} - V_p' \cdot \frac{d}{du}$
is stretched exponential decaying to equilibrium. More precisely,
$$
\var_\nu( \mathbf{S}_t f ) \leq \frac{1}{c} e^{-c t^{p/(2-p)}} \mathrm{Osc}(f)^2  \qquad \forall f, \forall t>0
$$
for some constant $c=c(p)$ (the lack of smoothness of $V_p$ at $0$ is just little nuisance that one can easily handle).

In \cite{bobkov-zeg2}, Bobkov and Zegarlinski proved, using weak Poincar\'e inequalities, 
that some Gibbs measures in infinite volume (with self potential $V_p$)
also satisfies a stretched exponential decay as above, but with 
a worst exponent. In this paper we shall prove the correct stretched exponential decay with the exponent 
$t^{p/(2-p)}$ not only for the potential $V_p$ but also for a larger class of potentials
of sub-exponential type. Moreover, our technique, based on the bisection approach \cite{martinelli}
together with the quasi factorisation property of the variance \cite{bcc}, applies also to potentials 
of the type $V=(1+\alpha)\log(1+|u|)$, $\alpha>0$ leading to Cauchy type distributions and polynomial decay to equilibrium.

Note that there is a difficulty here with respect to the usual Poincar\'e and/or log Sobolev inequalities.
Namely, weak Poincar\'e inequalities do not tensorise in general. In turn, there is no hope for
a dimension free analysis, and one has to take care of the growing dimension 
(see Section \ref{sec:wp} for a discussion about this fact).

\bigskip

The paper is organised as follows.
The notations and the setting, in particular the Hamiltonian and the Gibbs measure we consider, are introduced in the next section. Section \ref{sec:wp} is dedicated to the weak Poincar\'e inequalities, we recall few known facts
and prove some perturbation properties.
The results about the infinite volume Gibbs measure are collected and proved in Section \ref{sec:results}.
The main ingredients used in the proof of our theorems are postponed for the clarity of the exposition
to the last two sections 

\bigskip

\subsection*{Acknowledgement}
We would like to thank Fabio Martinelli, Senya Shlosman, Nobuo Yoshida, Boguslaw Zegarlinski and Pierre-Andr\'e Zitt for
some usefull discussions on the topic of this work.

\section{Notations and setting.}

\subsection{The Configuration space.}
The \emph{configuration space} we consider is $\Omega=\mathbb{R}^{\mathbb{Z}^d}$
where $d \geq 1$ is an integer that denotes the dimension of the lattice $\mathbb{Z}^d$.
Given $\Lambda \Subset \mathbb{Z}^d$ ({\it i.e.} $\Lambda$ is a finite subset of $\mathbb{Z}^d$),
we shall also deal with $\Omega_\Lambda = \mathbb{R}^\Lambda$. For any \emph{configuration} $\sigma \in \Omega$,
any site $x \in \mathbb{Z}^d$ and any $\Lambda \Subset \mathbb{Z}^d$, $\sigma_x$ stands for the value of the configuration (or the spin) at $x$ while $\sigma_\Lambda$ is the configuration $\sigma$ restricted to $\Lambda$.
We denote by $\mathcal{B}_\Lambda$ the $\sigma$-algebra of all Borell sets of $\Omega_\Lambda$.

A function which is measurable with respect to $\mathcal{B}_\Lambda$ with $\Lambda \Subset \mathbb{Z}^d$ is said to be \emph{local}. For any smooth local function, we set $|||f||| = \sum_{x \in \mathbb{Z}^d} \| \nabla_x f \|_\infty$
where $\| \cdot \|_\infty$ is the sup norm and $\nabla_x$ denotes the 
derivative with respect to the variable $\sigma_x$.

The Euclidean distance on $\mathbb{Z}^d$ is denoted by $d$.
With a slight abuse of notation, for $a,b \in \mathbb{Z}$, we shall often set $[ a,b ]$ for $[a,b]\cap \mathbb{Z}$.

\subsection{The Hamiltonian and the potentials.} \label{sec:h}
For any  $\Lambda \Subset \mathbb{Z}^d$  the \emph{Hamiltonian} $H_\Lambda : \Omega \to \mathbb{R}$ is given by
\begin{equation} \label{eq:hamiltonian}
H_\Lambda ( \sigma ) = 
\sum_{x \in \Lambda}  V(\sigma_x) + 
\frac{1}{T}  \!\!\! \sum_{y : d(x,y) \leq r} \!\!\! W(\sigma_x - \sigma_y) 
\qquad \forall \sigma \in \Omega 
\end{equation}
where $V, W : \mathbb{R} \to \mathbb{R}$ correspond respectively to the \emph{self potential} and the \emph{interaction potential} (pair potential).
The parameter $T \in (0,\infty)$ is the \emph{the temperature} 
and $r \in \mathbb{N}\setminus \{0\}$ is the \emph{range} of the interaction.
For $\sigma, \tau \in \Omega$, we also let $H_\Lambda^\tau(\sigma) := H_\Lambda (\sigma_\Lambda \tau_{\Lambda^c})$,
where $\sigma_\Lambda \tau_{\Lambda^c}$ stands for the configuration equal to $\sigma$ on $\Lambda$ and to $\tau$ on $\Lambda^c$ (the complement of $\Lambda$), and $\tau$ is called the \emph{boundary condition}.

\begin{remark}
One could consider more general Hamiltonian $H_\Lambda$ with {\it e.g.} infinite range interactions and/or
interaction potentials $W$ depending on the values of more than only two spins and/or depending on the sites etc. All the results below would hold in those more general settings, under specific assumptions. 
We make the choice of dealing with the Hamiltonian \eqref{eq:hamiltonian} for simplicity and for the clarity of the exposition.
\end{remark}

\smallskip
Now we describe our assumptions on $V$ and $W$.
We collect in Hypothesis $(H1)$ some smoothness conditions on $V$ and $W$.
\begin{assumption}[H1]
Given a self potential $V: \mathbb{R} \to \mathbb{R}$ and an interaction potential $W : \mathbb{R} \to \mathbb{R}$,
we say that Hypothesis $(H1)$ is satisfied if
\begin{itemize}
 \item $V$ is $\mathcal{C}^1$ and $\int_\mathbb{R} e^{-V(u)}du < \infty$; 
 \item $W$ is twice differentiable, $\|W\|_\infty<\infty$, $\|W'\|_\infty < \infty$ and $\|W''\|_\infty < \infty$.
\end{itemize}
\end{assumption}
The second assumption on $V$ guarantees that $d\nu(u)=Z_V^{-1} e^{-V}$ (with $Z_V=\int e^{-V}$) defines a probability measure.
The smoothness assumption about $V$ will be needed when defining the Glauber dynamics.

On the other side, the assumptions about $W$ will be usefull when defining the infinite volume Gibbs measure.

More specifically, the self potentials $V : \mathbb{R} \to \mathbb{R}$ we shall consider enter
into two classes of examples: $\kappa$-concave probability measures (a notion introduced by Borell \cite{borell75}, see
\cite{bob07} for a comprehensive introduction) and
sub-exponential like laws. More precisely, given any convex function $U : \mathbb{R} \to (0,\infty)$,
we shall consider either $V=(1+\alpha) \log U$  with $\alpha>0$ ($\kappa$-concave case with $\kappa= - 1/\alpha$)
or $V=U^p$, $p>0$ (sub-exponential like laws). 

The corresponding probability measure
 $d\nu(u) = Z_V^{-1} e^{-V(u)}du$ on $\mathbb{R}$ (with $Z_V = \int e^{-V(u)} du$) reads respectively as
$$
d\nu(u) = \frac{1}{Z_\alpha U(u)^{1+\alpha}} du 
\qquad \mbox{and} \qquad
d\nu(u) = Z_p^{-1}e^{-U(u)^{p}}du .
$$
Prototypes are respectively the Cauchy distributions ($U(u) = 1+|u|$ or equivalently 
$V(u)=(1+\alpha) \log(1+|u|)$) and  the sub-exponential laws ($U(u)=|u|$ or equivalently $V(u)=|u|^p$):
\begin{equation}\label{eqcauchy}
d\nu(u) = \frac{\alpha}{2\left(1+|u|\right)^{1+\alpha}}du
\qquad \mbox{and} \qquad
d\nu(u) = \frac{e^{-|u|^{p}}}{2\Gamma(1+\frac{1}{p})}du
\end{equation}
for $\alpha > 0$ and $p \in (0,1]$.
In both examples the measure $\nu$ has tails larger than exponential.
For the sub-exponential law, in order to fulfil Hypothesis $(H1)$ one could consider 
{\it e.g.} $U(x)=\sqrt{1+u^2}$ to avoid differentiability trouble in $0$.

\subsection{The Gibbs measures}
The \emph{finite volume Gibbs measure} in $\Lambda \Subset \mathbb{Z}^d$ at temperature $T$ 
and boundary condition $\tau$ is given by
\begin{equation} \label{eq:gibbs}
\mu_\Lambda^\tau (d\sigma) = \left( Z_\Lambda^\tau \right)^{-1} \exp \left\{ -H_\Lambda^\tau(\sigma) \right\} \prod_{x \in \Lambda} d\sigma_x \times 
\delta_{\Lambda^c, \tau}( d\sigma)
\end{equation}
where $\delta_{\Lambda^c, \tau}$ is the Dirac probability measure on $\Omega_{\Lambda^c}$ which gives mass $1$ to the configuration $\tau$ and $Z_\Lambda^\tau$ is the proper normalisation factor.
We denote with $\mu_\Lambda^\tau(f)$ the expectation of $f$ with respect to $\mu_\Lambda^\tau$, while 
$\mu_\Lambda(f)$ denotes the functions $\tau \mapsto \mu_\Lambda^\tau (f)$. For any Borel set $X \subset \Omega$
we set $\mu_\Lambda(X):=\mu_\Lambda( \mathds{1}_X)$, where $\mathds{1}_X$ is the characteristic function on $X$.
We write
$\mu_\Lambda(f,g)$ to denote the covariance (with respect to $\mu_\Lambda$) of $f$ and $g$
and $\var_{\mu_\Lambda}(f)=\mu_\Lambda(f,f)$ for the variance of $f$ under $\mu_\Lambda$.

The family of measures \eqref{eq:gibbs} satisfies the DLR compatibility conditions: for all 
Borell sets $X \subset \Omega$
$$
\mu_\Delta (\mu_\Lambda (X)) = \mu_\Delta (X) \qquad \qquad \forall \Lambda, \Delta \Subset \mathbb{Z}^d \mbox{ such that }
\Lambda \subset \Delta. 
$$
If in addition of Hypothesis $(H1)$, $T$ is large enough, then (see \cite[Proposition (8.8)]{georgii}) the Dobrushin's uniqueness condition is satisfied. Hence there exists a unique
infinite volume Gibbs measure $\mu$ satisfying $\mu(\mu_\Lambda (X))=\mu(X)$ for $\Lambda \Subset \mathbb{Z}^d$ and
any Borell set $X \subset \Omega$.
Moreover (see \cite[Remark (8.26) together with Corollary (8.32)]{georgii}) (if $T$ is large enough) there exist constants $D=D(r,T,\|W\|_\infty)$ and $m=m(r,T,\|W\|_\infty)$ 
such that for any $\Lambda \Subset \mathbb{Z}^d$, it holds 
\begin{equation} \label{eq:smc}
|\mu_\Lambda^\tau (f,g)| \leq D |\Delta_f| |\Delta_g| \|f \|_\infty \|g\|_\infty e^{-md(\Delta_f , \Delta_g)} 
\end{equation}
for any boundary condition $\tau$, any bounded local functions $f$ and $g$ with support $\Delta_f$ and $\Delta_g$ satisfying  $\Delta_f, \Delta_g \subset \Lambda$. Here $|\cdot|$ stands for the Lebesgue measure on $\mathbb{Z}^d$. 
Inequality \eqref{eq:smc} is known as the  \emph{strong mixing condition}. 
Note that the previous argument does not depend on the self potential $V$. 

In the sequel, we will always assume the following: 
\begin{assumption}[H2]
Given the potentials $V$ and $W$, the temperature $T$, 
we say that Hypothesis $(H2)$ is satisfied if there exists a unique infinite volume Gibbs measure $\mu$ and if
the strong mixing condition \eqref{eq:smc} holds true.
\end{assumption}

In  particular, by the argument above, if $(H1)$ is satisfied then $(H2)$ is also satisfied as soon as $T$ is 
large enough, whatever the choice of $V$.

\subsection{The dynamics.}

The dynamics we consider are of Glauber type. For any $\Lambda \Subset \mathbb{Z}^d$,
any boundary condition $\tau \in \Omega$, let $(\PT{t}^{\Lambda,\tau})_{t \geq 0}$ be the Markov semi-group
associated to the generator
\begin{equation} \label{eq:gi}
 \GI_\Lambda^\tau = \sum_{x \in \Lambda} \Delta_x - \sum_{x \in \Lambda} \nabla_x H_\Lambda^\tau \cdot \nabla_x
\end{equation}
where %$\mathcal{L}_x = \Delta_x - V'(\sigma_x) \cdot \nabla_x$ is a one dimentional diffusion operator 
%acting only on the coordinate $\sigma_x$. Here 
$\nabla_x$ and $\Delta_x$ stand respectively for the first and second 
partial derivative with respect to the variable $\sigma_x$. When there is no confusion, we shall drop the superscript $\tau$ in the definition and write simply $\PT{t}^\Lambda$ and $\GI_\Lambda$. The generator 
$\GI_\Lambda^\tau$ is symmetric in $\mathbb{L}^2( \mu_\Lambda^\tau)$.
On the other hand, $\PT{t}^{\Lambda,\tau}$ is a contraction on 
$\mathbb{L}^p(\mu_\Lambda^\tau)$ for all $p\in[1,\infty]$.

For any $\Lambda \Subset \mathbb{Z}^d$ we denote by $\mathcal{D}_\Lambda^\tau$ 
(we shall also often drop the superscript $\tau$) the \emph{Dirichlet form} associated 
to $\GI_\Lambda^\tau$ and defined by
\begin{equation} \label{eq:dirichlet}
 \mathcal{D}_\Lambda^\tau (f) = \frac{1}{2} \sum_{x \in \Lambda} \mu_\Lambda^\tau \left( |\nabla_x f|^2 \right)
\end{equation}
for all sufficiently smooth $f$.

If $(H1)$ and $(H2)$ are satisfied,
it is possible (see {\it e.g.} \cite{guionnet-zegarlinski}) to construct the infinite volume semigroup
by
\begin{equation} \label{eq:semi-group}
\PT{t} = \lim_{\Lambda \to \mathbb{Z}^d} \PT{t}^\Lambda
\end{equation}
on the space of bounded  smooth functions (in particular $(H1)$ and $(H2)$ guarantee that the limit above does not depend on the boundary condition). The associated infinite volume generator will be denoted by $\GI$.

\section{The weak Poincar\'e inequalities.} \label{sec:wp}

The main tool we shall use is the so called weak Poincar\'e inequality. We now introduce this notion,
first on $\mathbb{R}$, and then on $\Omega_\Lambda$. We then collect some usefull results we shall use later.

\subsection{Introduction.}

In this section we introduce the notion of weak Poincar\'e inequality on $\mathbb{R}$. Then we derive some (known)
bounds on the decay to equilibrium. Finally we explain how to get weak Poincar\'e inequalities on product spaces.

\smallskip

Consider the probability measure $d\nu(u)=Z_V^{-1} e^{-V(u)}du$, on $\mathbb{R}$. We say that
$\nu$ satisfies a \emph{weak Poincar\'e inequality} with rate function $\beta : (0,\infty) \to [0,\infty)$,
if for any bounded function $f : \mathbb{R} \to \mathbb{R}$ smooth enough, it holds
\begin{equation} \label{eq:wp}
 \var_\nu(f) \leq \beta(s) \int (f')^2 d\nu + s \mathrm{Osc}(f)^2 \qquad \forall s >0 .
\end{equation}
where $\mathrm{Osc}(f)$ is \emph{the oscillation} of $f$: $\mathrm{Osc}(f):=\sup f - \inf f$. 

Note that if $\lim_{s \to 0} \beta(s) = \beta_0 < \infty$, then
Inequality \eqref{eq:wp} reduces to the standard Poincar\'e inequality
$$
\var_\nu (f) \leq \beta_0 \int (f')^2 d\nu .
$$
On the other hand, since $\var_\nu (f) \leq \frac{1}{4} \mathrm{Osc}(f)^2$, only the values $s \in (0,1/4)$ are relevant. 
Most of the information is encoded in the behaviour of  $\beta$ near the origin.
Weak Poincar\'e inequalities have been introduced by R\"ockner and Wang \cite{rockner-wang}.
One interested feature of Inequality \eqref{eq:wp} is that it gives a control on
the $\mathbb{L}^2$ decay to equilibrium of the Markov
semi-group $(\mathbf{S}_{t})_{t\geq 0}$ on $\mathbb{R}$ with generator $L=\frac{d^2}{du^2} - V' \cdot \frac{d}{du}$.
\begin{proposition}[\cite{rockner-wang}]\label{prop:rwang}
   Let $\nu$ be a probability measure on $\mathbb{R}$ with density $Z_V^{-1}e^{-V}$
   with respect to the Lebesgue measure on $\mathbb{R}$. 
   Let $(\mathbf{S}_{t})_{t \ge 0}$ be the corresponding semi-group
   with generator $L:= \frac{d^2}{du^2} - V' \cdot \frac{d}{du}$.
   If $\nu$ satisfies the weak Poincar\'e inequality \eqref{eq:wp} with rate function $\beta$, then,
   every smooth $f : \mathbb{R} \rightarrow \mathbb{R}$ satisfies
   \begin{equation} \label{eq:decay}
   \var_\nu (\mathbf{S}_{t}f) \leq e^{-\frac{2t}{\beta(s)}} \var_\nu (f)
   +4s(1-e^{-\frac{2t}{\beta(s)}}) \|f - \nu(f)\|_\infty^2   \qquad \forall s,t >0 .
   \end{equation}
\end{proposition}

The result by R\"ockner and Wang holds in more general settings, see  \cite{rockner-wang}.
We sketch the proof for completeness.
\begin{proof}
   Assume without loss of generality that $\nu(f)=0$ (which implies $\nu(\mathbf{S}_{t}f)=0$ for any $t$).
   If $u(t)= \var_\nu( \mathbf{S}_{t}f) = \int (\mathbf{S}_{t}f)^2 d\nu$, the weak Poincar\'e inequality implies
   that
   $$
   u'(t) = 2 \int  \mathbf{S}_{t}f L \PT{t}f \, d\nu=- 2 \int |\frac{d}{du} \mathbf{S}_{t}f|^{2}d\nu \leq
   -\frac{2}{\beta(s)}\left[u(t) - 2s \|f\|_\infty\right]
   $$
   since $\mathrm{Osc}(\mathbf{S}_{t} f) \leq 2 \| \mathbf{S}_{t}f\|_\infty \leq 2\|f\|_\infty$. The result
   follows by integration.
\end{proof}

For the two classes of self-potentials $V$ introduced above, the corresponding rate function $\beta$ has been computed 
in \cite{cggr} (see also \cite{rockner-wang,bcr2,bobkov-zeg2,bob07}). Given $\beta$, one can 
then optimise over $s>0$ in \eqref{eq:decay} to get an explicit decay of the Markov semi-group 
$(\mathbf{S}_{t})_{t \geq 0}$
in $\mathbb{L}^2(\nu)$.
Let $U : \mathbb{R} \to (0,\infty)$ be a convex function.
\begin{itemize}
 \item 
If $V=(1+\alpha) \log U$  with $\alpha>0$ (the $\kappa$-concave case), then $\nu$
satisfies a weak Poincar\'e inequality with rate function $\beta(s) = c_\alpha s^{-2/\alpha}$ for some constant 
$c_\alpha>0$ (see  \cite[Proposition 5.4]{cggr}).
Optimising \eqref{eq:decay} over $s$ (together with some computations given in \cite[Corollary 2.4]{rockner-wang},
see also the proof of Corollary \ref{cor:kconcave} below) 
leads to
\begin{equation} \label{eq:decay1}
\var_\nu(\mathbf{S}_{t}f) \leq \frac{C}{t^{\alpha/2}} \| f - \nu(f) \|_\infty^2
\end{equation}
for some constant $C=C(\alpha)>0$.

\item If $V=U^p$, $p \in (0,1)$ (the sub-exponential case), then $\nu$
satisfies a weak Poincar\'e inequality with rate function 
$\beta(s) = c_p \left( \log \frac{2}{s \wedge 1} \right)^{2(\frac{1}{p}-1)}$ for some constant $c_p>0$ (see 
\cite[Proposition 5.6]{cggr}). 
Optimising \eqref{eq:decay} over $s$ (take $s=e^{-ct^{p/(2-p)}}$) leads to
\begin{equation} \label{eq:decay2}
\var_\nu(\mathbf{S}_{t}f) \leq \frac{1}{C} e^{-Ct^{p/(2-p)}} \| f - \nu(f) \|_\infty^2
\end{equation}
for some constants $C=C(p)>0$. Note that $p/(2-p) \in (0,1)$.
\end{itemize}

The previous results are optimal, in the sens that for $U(u)=1+|u|$, respectively $U(u)=|u|$, 
neither the rate function $\beta$ nor the $\mathbb{L}^2$ decay can be improved.
In particular there is no hope for a Poincar\'e inequality to hold, or equivalently, for an exponential decay
to equilibrium in $\mathbb{L}^2$.

Note that the limiting case $p=1$ corresponds to the exponential measure for which it is known that a Poincar\'e
inequality holds, and thus an exponential decay of the semi-group. This fact is encoded in the rate function $\beta$
(which becomes a constant) and on the decay \eqref{eq:decay2} (which becomes exponential).

\bigskip

Contrary to the Poincar\'e inequality, the weak Poincar\'e inequalities do not tensorise in general. If the probability measure
$\nu_n=\otimes\nu^{(i)}$ on $\mathbb{R}^n$ is the tensor product of $n$ copies of $\nu$, it is possible 
(and actually very easy, see \cite[Section 3]{bcr2}) to prove that
\begin{equation} \label{eq:product}
\var_{\nu_n}(f) \leq \beta(s/n) \int \sum_{i=1}^n |\nabla_i f|^2 d\nu_n + s \mathrm{Osc}(f)^2
\qquad \forall s>0
\end{equation}
for all functions $f : \mathbb{R}^n \to \mathbb{R}$ smooth enough. The rate function $\beta(\cdot/n)$ is best possible
for the product of Cauchy measures and sub-exponential laws introduced in \eqref{eqcauchy}.
In particular, there is no hope for those measures with heavy tails to get a weak Poincar\'e inequality in infinite dimension. A deep explanation of this phenomenon can be found in
Talagrand's paper \cite{Tal91} (see also the introduction of \cite{bcr2}). It relies on the 
concentration of measure phenomenon.

However, quite remarkable is the fact that the decays \eqref{eq:decay1} and \eqref{eq:decay2} 
still hold in the infinite system $\Omega$ with 
infinite volume Gibbs measure $\mu$ and Markov semi-group $(\PT{t})_{t \geq 0}$ introduce in the previous section. 
The aim of this paper is to prove such results.

\subsection{The weak Poincar\'e inequalities for Gibbs measures.}

Now we turn to the Gibbs measure setting of the previous section. 
Let $\Lambda \Subset \mathbb{Z}^d$.
For any $s>0$, let $\beta_\Lambda(s)$ 
be the best non-negative number such that for any boundary condition $\tau$ and any smooth function $f : \Omega_\Lambda \to \mathbb{R}$,
\begin{equation} \label{eq:defb}
 \var_{\mu_\Lambda^\tau}(f) \leq \beta_\Lambda(s) \mathcal{D}_\Lambda^\tau (f) + s \mathrm{Osc}(f)^2 .
\end{equation}
By this procedure we have defined a non-increasing function $\beta_\Lambda : (0,\infty) \to [0,\infty)$. 
Note that the system is invariant under translation and rotation. 
Hence, two finite subsets of $\mathbb{Z}^d$ that are equal
under translation and rotation lead to the same rate function $\beta$.

\smallskip

Our aim is to get the best possible rate function for \eqref{eq:defb} to hold. In view of \eqref{eq:product},
the best one can hope is $\beta(s/|\Lambda|)$ if $\beta$ denotes the rate function associated to
the one dimensional measure $d\nu=Z_V^{-1}e^{-V}$. This will actually be almost true, see 
Proposition \ref{prop:perturbation2} below. The difficulty here comes from the interacting part which 
can be of order $e^{|\Lambda|}$. The following perturbation result goes in this direction.
Even if it is far from being optimal, it will be usefull in the proof of our main result.

\begin{proposition}[Perturbation] \label{prop:perturbation}
Assume $(H1)$.
Also, assume that the self-potential $V$ is such that
$d\nu=Z_V^{-1}e^{-V}$ satisfies the following 
weak Poincar\'e inequality on $\mathbb{R}$ for some non-increasing rate function $\beta$:
$$
\var_\nu(f) \leq \beta(s) \int (f')^2 d\nu + s \mathrm{Osc}(f)^2 \qquad \forall f, \;\forall s >0 .
$$
Then, there exists a constant $C=C(r,T,d,\|W\|_\infty)$ such that 
for any $\Lambda \Subset \mathbb{Z}^d$, any boundary condition $\tau \in \Omega$, any smooth $f : \Omega_\Lambda \to \mathbb{R}$,
$$
\var_{\mu_\Lambda^\tau}(f) \leq C e^{C|\Lambda|} \beta\left(\frac{s}{C|\Lambda|e^{C|\Lambda|}}\right) \mathcal{D}_\Lambda^\tau (f) + s \mathrm{Osc}(f)^2 
\qquad \forall s>0.
$$
\end{proposition}

\begin{remark}
A somehow similar statement can be found in \cite[Lemma 12.1]{bobkov-zeg2}.
\end{remark}

\begin{proof}
Fix $\Lambda \Subset \mathbb{Z}^d$ and $\tau \in \Omega$.
Let $d\nu_\Lambda (\sigma) = Z_V^{-|\Lambda|} \exp\{ - \sum_{x \in \Lambda} V(\sigma_x) \}d\sigma_\Lambda$ be the 
probability measure corresponding to the product part of $\mu_\Lambda^\tau$.
By the product property \eqref{eq:product} we have
$$
\var_{\nu_\Lambda}(f) \leq \beta\left(\frac{s}{|\Lambda|}\right) \sum_{x\in\Lambda} \nu_\Lambda \left( |\nabla_x f|^2 \right) 
+ s \mathrm{Osc}(f)^2 \qquad \forall f:\Omega_\Lambda \to \mathbb{R}, \; \forall s>0 .
$$
Now Hypothesis $(H1)$ guarantees that there exists a constant $C$ (depending on $r$, $T$, $d$ and  $\|W\|_\infty$ but independent of the boundary condition $\tau$ and $\Lambda$)
such that
$$
C^{-1} e^{-C|\Lambda|}
\leq 
\frac{\nu_\Lambda(\sigma)}{\mu_\Lambda^\tau(\sigma)}
\leq C e^{C|\Lambda|}
\qquad \forall \sigma \in \Omega_\Lambda .
$$
Hence, since  $\var_{\mu_\Lambda^\tau}(f)=\inf_a \mu_\Lambda^\tau ((f-a)^2)$, we get for any $s>0$,
\begin{eqnarray*}
\var_{\mu_\Lambda^\tau}(f) 
& \leq & 
C e^{C|\Lambda|} \var_{\nu_\Lambda}(f) \\
& \leq &
C e^{C|\Lambda|} \beta \left(\frac{s}{|\Lambda|}\right) \sum_{x\in\Lambda} \nu_\Lambda \left( |\nabla_x f|^2 \right) 
+ s C e^{C|\Lambda|} \mathrm{Osc}(f)^2 \\
& \leq &
2 C^2 e^{2C|\Lambda|} \beta \left(\frac{s}{|\Lambda|}\right) \mathcal{D}_\Lambda^\tau (f)
+ s C e^{C|\Lambda|} \mathrm{Osc}(f)^2 .
\end{eqnarray*}
The result follows.
\end{proof}

As for the tensorisation property, the previous result is of no help in order to get directly 
infinite volume estimates, since when $|\Lambda| \to \infty$ the weak Poincar\'e inequality becomes trivial.

Using the bisection technique \cite{martinelli}, the result of Proposition \ref{prop:perturbation} can be improved
for volumes $\Lambda$ that are cubes. More precisely we have the following proposition.

\begin{proposition}[Perturbation improved] \label{prop:perturbation2}
Assume $(H1)$ and $(H2)$. 
Also, assume that the self-potential $V$ is such that
$d\nu=Z_V^{-1}e^{-V}$ satisfies the following 
weak Poincar\'e inequality on $\mathbb{R}$ for some non-increasing rate function $\beta$:
$$
\var_\nu(f) \leq \beta(s) \int (f')^2 d\nu + s \mathrm{Osc}(f)^2 \qquad \forall f, \;\forall s >0 .
$$
Then, for any $\varepsilon \in (0,1)$, there exists a constant $C=C(\varepsilon,r,T,d,\|W\|_\infty)$
such that for any integer $L$,
\begin{equation} \label{eq:pert}
\var_{\mu_\Lambda^\tau} (f) \leq C \beta \left( \frac{s}{C|\Lambda|^{1 + \varepsilon}} \right) \mathcal{D}_\Lambda^\tau (f) 
+ s \mathrm{Osc}(f) \qquad \forall \tau \in \Omega, \; \forall f, \; \forall s>0,
\end{equation}
where $\Lambda=[-L,L]^d$.
\end{proposition}
The proof is postponed to Section \ref{sec:improved}.
Note that we obtained a quasi optimal inequality, up to the power $\varepsilon$. Indeed Inequality \eqref{eq:pert}
is closed to the non-interacting case \eqref{eq:product}.
Proposition \ref{prop:perturbation2} is at the heart of
the proof of the main theorems through the following two Lemmas.
In fact, using the semi-group strategy \cite{rockner-wang} explained in Proposition \ref{prop:rwang}, 
Inequality \eqref{eq:pert} already leads to some
finite volume decay of the semi-group $(\PT{t}^{\Lambda,\tau})_{t \geq 0}$, for cubes.

\begin{corollary}[$\kappa$-concave case] \label{cor:kconcave}
Let $U : \mathbb{R} \to (0,\infty)$ be a convex function and  $V=(1+\alpha) \log U$  with $\alpha>0$.
Assume $(H1)$ and $(H2)$. Then, for any $\varepsilon$, there exists a constant $C=C(\varepsilon,\alpha,r,T,d,\|W\|_\infty)$
such that for any integer $L$, any local function $f$ satisfies
$$
\var_{\mu_\Lambda^\tau} (\PT{t}^{\Lambda,\tau}f) 
\leq 
C \frac{|\Lambda|^{1+\varepsilon}}{t^{\alpha/2}} \|f - \mu_\Lambda^\tau(f) \|_\infty^2
\qquad \forall t>0, \; \forall \tau \in \Omega,
$$
where $\Lambda=[-L,L]^d$. 
\end{corollary}

\begin{proof}
Fix an integer $L$, $\tau \in \Omega$, $\varepsilon >0$ and a local function $f$. Set $\Lambda=[-L,L]^d$.
Assume without loss of generality that $\mu_\Lambda^\tau(f)=0$.

As mentioned before Inequality \eqref{eq:decay1}, the measure 
$d\nu=Z_V^{-1}e^{-V}$ on $\mathbb{R}$ satisfies a weak Poincar\'e inequality with rate function
$\beta(s) = c_\alpha s^{-2/\alpha}$ for some constant $c_\alpha>0$.
Hence, using Proposition \ref{prop:perturbation2}, $\mu_\Lambda^\tau$ satisfies a weak Poincar\'e 
inequality with rate function $\gamma(s)= C s^{-2/\alpha} |\Lambda|^{2(1+\varepsilon)/\alpha}$, for some constant 
$C=C(\varepsilon,\alpha,r,T,d,\|W\|_\infty)$. In turn, using the strategy of the proof of Proposition 
\ref{prop:rwang} (we omit the details), we get that 
$$
\var_{\mu_\Lambda^\tau} (\PT{t}^{\Lambda,\tau}f) \leq 
e^{-\frac{2t}{\gamma(s)}} \var_\nu (f)
   +4s(1-e^{-\frac{2t}{\gamma(s)}}) \|f \|_\infty^2   \qquad \forall s,t >0 .
$$
Following \cite{rockner-wang}, we take $s=(\lambda/t)^{\alpha/2}$ with $\lambda>0$ chosen in such a way that
$$
e^{-\frac{2t}{\gamma(s)}} = e^{-\frac{2 t s^{2/\alpha}}{C |\Lambda|^{2(1+\varepsilon)/\alpha}}} = e^{-\frac{2 \lambda}{C |\Lambda|^{2(1+\varepsilon)/\alpha}}} = \left(\frac{1}{2} \right)^{\frac{\alpha}{2}+1} .
$$
It follows that
$$
\var_{\mu_\Lambda^\tau} (\PT{t}^{\Lambda,\tau}f) \leq \left( \frac{1}{2} \right)^{\frac{\alpha}{2}+1} \var_\nu (f) 
+ 4 \left( \frac{\lambda}{t} \right)^{\frac{\alpha}{2}} \|f  \|_\infty^2 \qquad \forall t>0 .
$$
We omit the superscript $\tau$. Applying this inequality repeatedly, we obtain 
(using also the fact that $\PT{t}^{\Lambda}$ is a contaction in 
the sup-norm)
\begin{align*}
\var_{\mu_\Lambda} (\PT{t}^{\Lambda}f) & = 
\var_{\mu_\Lambda} \left(\PT{t/2}^{\Lambda} \left[\PT{t/2}^{\Lambda}f \right] \right) \\
& \leq 
\left( \frac{1}{2} \right)^{\frac{\alpha}{2}+1} \var_\nu \left(\PT{t/2}^{\Lambda}f \right) 
+ 4 \left( \frac{\lambda}{t} \right)^{\frac{\alpha}{2}} 2^\frac{\alpha}{2} \|f  \|_\infty^2 \\
& \leq
\left( \frac{1}{4} \right)^{\frac{\alpha}{2}+1} \var_\nu \left(\PT{t/4}^{\Lambda}f \right) 
+ 4 \left( \frac{\lambda}{t} \right)^{\frac{\alpha}{2}} 2^\frac{\alpha}{2} \|f  \|_\infty^2 \left(1+\frac{1}{2}\right) \\
& \leq \cdots \\
& \leq 
4 \left( \frac{\lambda}{t} \right)^{\frac{\alpha}{2}} 2^\frac{\alpha}{2} \|f  \|_\infty^2
\sum_{n \geq 0} 2^{-n} .
\end{align*}
The result follows by our choice of $\lambda$.
\end{proof}

The next result concerns sub-exponential type laws.

\begin{corollary}[Sub-exponential case] \label{cor:subex}
Let $U : \mathbb{R} \to (0,\infty)$ be a convex function and  $V= |U|^p$  with $p \in (0,1)$.
Assume $(H1)$ and $(H2)$. Fix $A >0$. Then,
there exists a constant $C=C(p,A,r,T,d,\|W\|_\infty)$
such that for any integer $L$, any local function $f$ satisfies
$$
\var_{\mu_\Lambda^\tau} (\PT{t}^{\Lambda,\tau}f) 
\leq 
\frac{1}{C} e^{-C t^{p/(2-p)}} \|f - \mu_\Lambda^\tau(f) \|_\infty^2
\qquad \forall \tau \in \Omega, 
$$
 provided $t^{p/(2-p)} \geq 2A \log(|\Lambda|)$, where $\Lambda=[-L,L]^d$.
\end{corollary}

\begin{proof}
 Fix an integer $L$, $\tau \in \Omega$ and a local function $f$ with $\mu_\Lambda^\tau(f)=0$.
We start as in the proof of Corollary \ref{cor:kconcave}, using instead that the one dimensional measure 
$d\nu=Z_V^{-1}e^{-V}$  satisfies a weak Poincar\'e inequality with rate function
$\beta(s) = c_p \left( \log \frac{2}{s \wedge 1} \right)^{2(\frac{1}{p}-1)}$ for some constant $c_p>0$
(see before Inequality \eqref{eq:decay2}), to get that
$$
\var_{\mu_\Lambda^\tau} (\PT{t}^{\Lambda,\tau}f) \leq 
e^{-\frac{2t}{\gamma(s)}} \var_\nu (f)
   +4s(1-e^{-\frac{2t}{\gamma(s)}}) \|f \|_\infty^2   \qquad \forall s,t >0 
$$
with $\gamma(s) = C \left( \log \frac{2|\Lambda|^{3/2}}{s \wedge |\Lambda|^{3/2}} \right)^{2(\frac{1}{p}-1)}$
for some $C=C(p,r,T,d,\|W\|_\infty)$ (we have chosen $\varepsilon=3/2$ in Proposition \ref{prop:perturbation2}).

Choose $s= e^{- t^{p/(2-p)}}$. Under the assumption $t^{p/(2-p)} \geq A \log(|\Lambda|^{3/2})$,
the expected result follows after few rearrangements.
\end{proof}

%%%%%%%%%%%%%%%%%%%%%%%%%%%%%%%%%%%%%%%%%%%%%%%%%%%%%%%%%%%%%%%%%%%%%%%%%%%%%%%%%%%%%%%%%%%%%%%%%%
%%%%%%%%%%%%%%%%%%%%%%%%%%%%%%%%%%%%%%%%%%%%%%%%%%%%%%%%%%%%%%%%%%%%%%%%%%%%%%%%%%%%%%%%%%%%%%%%%%
%%%%%%%%%%%%%%%%%%%%%%%%%%%%%%%%%%%%%%%%%%%%%%%%%%%%%%%%%%%%%%%%%%%%%%%%%%%%%%%%%%%%%%%%%%%%%%%%%%

\section{The results.} \label{sec:results}

In this section we shall deal with the two classes of examples of self potential $V$ introduced in Section \ref{sec:h}. 
We assume that Hypothesis $(H1)$ and $(H2)$ are satisfied in such a way that the infinite
volume Gibbs measure $\mu$ exists. Recall that the Markov semi-group $(\PT{t})_{t \geq 0}$ has been 
defined in \eqref{eq:semi-group}. Also, we set $\| f \| = |||f||| + \|f\|_\infty$.

\begin{theorem}[$\kappa$-concave case] \label{th:kconcave}
Let $U : \mathbb{R} \to (0,\infty)$ be a convex function and $\alpha >0$. 
Set $V=(1+\alpha) \log U$. Let $W : \mathbb{R} \to \mathbb{R}$. Assume $(H1)$ and $(H2)$.
Fix and integer $\ell \geq 1$.
%Define the family of finite volume Gibbs measure by \eqref{eq:gibbs}. Assume that $(H2)$ is satisfied and denote
%by $\mu$ the infinite volume Gibbs measure on $\Omega$. Then, the Markov semi-group $(\PT{t})_{t \geq 0}$
%defined in \eqref{eq:semi-group} satisfies 
Then, for any $\varepsilon \in (0,1)$, there exists a constant 
$C$ depending on $\varepsilon$, $\ell$, $\alpha$, $T$, $r$, $d$, $\|W\|_\infty$, $\|W'\|_\infty$ and $\|W''\|_\infty$ 
such that for all bounded local functions $f : \Omega \to \mathbb{R}$ with $|\Delta_f| \leq \ell^d$,
\begin{equation} \label{eq:thkconcave}
 \var_\mu(\PT{t}f) \leq \frac{C}{t^{\frac{\alpha}{2} -d(1+ \varepsilon)}} \| f - \mu(f) \|^2 
\qquad \forall t>0.
\end{equation}
\end{theorem}

\begin{remark}
The spurious term $d(1+ \varepsilon)$ is, a priori, technical and we believe that the correct decay should be
with the exponent $\alpha/2$ as in the one dimensional case. But on the other hand, it could be that the very heavy tails
of the Cauchy type distributions slow down the dynamics with some strange unattended phenomenon
(that we  have not been able to catch).

Observe also that we obtain a polynomial decay only for  $\alpha>2d$. For $\alpha \leq 2d$, the previous bound is useless
since $\PT{t}$ is a contaction: we already know that $\var_\mu(\PT{t}f) \leq \| f - \mu(f) \|_\infty^2$.
\end{remark}

Similarly, we have for sub-exponential self-potentials:

\begin{theorem}[Sub-exponential case] \label{th:subexp}
Let $U : \mathbb{R} \to (0,\infty)$ be a convex function and $p \in (0,1)$. 
Set $V= |U|^p$. Let $W : \mathbb{R} \to \mathbb{R}$. Assume $(H1)$ and $(H2)$.
Fix an integer $\ell \geq1$.
%Define the family of finite volume Gibbs measure by \eqref{eq:gibbs}. Assume that $(H2)$ is satisfied and denote
%by $\mu$ the infinite volume Gibbs measure on $\Omega$. Then, the Markov semi-group $(\PT{t})_{t \geq 0}$
%defined in \eqref{eq:semi-group} satisfies 
Then, there exists a constant 
$C$ depending on $p$, $\ell$, $T$, $r$, $d$, $\|W\|_\infty$, $\|W'\|_\infty$ and $\|W''\|_\infty$ 
such that for all bounded local functions $f : \Omega \to \mathbb{R}$ with $|\Delta_f| \leq \ell^d$,
\begin{equation} \label{eq:thsubexp}
 \var_\mu(\PT{t}f) \leq \frac{1}{C} e^{-C t^{p/(2-p)}} \| f - \mu(f) \|^2   \qquad \forall t>0.
\end{equation}
\end{theorem}

The proof of Theorem \ref{th:kconcave} and Theorem \ref{th:subexp} relies on two main ingredients: 
the bisection technique \cite{martinelli} through Corollary \ref{cor:kconcave} and Corollary
\ref{cor:subex}, and the following property known as finite speed of propagation.

\begin{proposition}[Finite speed of propagation] \label{prop:finitespeed}
 Fix and integer $\ell \geq 1$ and assume $(H_1)$ and $(H2)$.
Then, for any local function $f$ with support $\Delta_f \subset [-\ell,\ell]^d$, 
any $L$ multiple of $r$, any boundary condition $\tau \in \Omega$,
$$
\| \PT{t}f - \PT{t}^{\Lambda,\tau} f \|_\infty \leq C ||| f ||| \left( \frac{C't}{L} \right)^{C''L} e^{Ct} 
\qquad \forall t>0
$$
with $\Lambda = [-L,L]^d$,
for some constant $C, C', C''>0$ depending only on $r$, $d$, $\|W'\|_\infty$, $\|W''\|_\infty$ and $\ell$.
\end{proposition}
The proof of Proposition \ref{prop:finitespeed} is postponed to Section \ref{sec:fs}
for the clarity of the exposition.

\begin{proof}[Proof of Theorem \ref{th:kconcave}]
Fix $t > 0$ and $\varepsilon \in (0,1)$.
Let $f$ be a local function. Since the system is invariant under translation and rotation we may assume
as we shall that the support $\Delta_f$ of $f$ contains the origin $0 \in \mathbb{Z}^d$. Furthermore, we can assume that 
$\mu(f)=0$.

Let $\Lambda = [-L,L]^d \Subset \mathbb{Z}^d$ with $L= \lambda t + \lambda'$, where $\lambda,\lambda'>0$ are  parameters
that will be chosen later. We assume that $\lambda'$ is large enough in such a way that $\Delta_f \subset \Lambda$.
Let $\tau \in \Omega$ be a boundary condition.
Our starting point is the following bound
\begin{equation} \label{eq:start}
\var_\mu (\PT{t}f) \leq 2 \| \PT{t}f - \PT{t}^{\Lambda,\tau} f \|_\infty^2 
+ 2 \var_{\mu_\Lambda^\tau} \left( \PT{t}^{\Lambda,\tau}f \right) .
\end{equation}
The first term of \eqref{eq:start} is controlled by the finite speed of propagation result above.
Indeed, we can choose $\lambda$ and $\lambda'$ large enough in such a way that 
$L$ is a multiple of $r$ and 
$$
\left( \frac{C't}{L} \right)^{C''L} e^{Ct} \leq e^{-ct}
$$
for some constant $c$ depending on $C$, $C'$, $C''$, $\lambda$ and $\lambda'$, where 
$C$, $C'$ and $C''$ are defined in Proposition \ref{prop:finitespeed}.

Going back to \eqref{eq:start} and thanks to Corollary \ref{cor:kconcave}, we get that
$$
\var_\mu (\PT{t}f) \leq  \frac{||| f |||}{c} e^{-ct} 
+ \frac{1}{c} \frac{|\Lambda|^{1+\varepsilon}}{t^{\alpha/2}} \| f \|_\infty^2 
$$
for some constant $c$ depending only on 
$\varepsilon$, $\alpha$, $T$, $r$, $d$, $\|W\|_\infty$, $\|W'\|_\infty$, $\|W''\|_\infty$ and $\ell$.
The expected result follows.
\end{proof}

\begin{proof}[Proof of Theorem \ref{th:subexp}]
The proof of Theorem \ref{th:subexp} is identical to the one of Theorem \ref{th:kconcave}. 
We use the same notations.

Note that, for $L=\lambda t + \lambda'$ and $\Lambda=[-L,L]^d$,
 there exists $A=A(\lambda,\lambda',p)$ such that for any $t \geq 1$,
$t^{p/(2-p)} \geq 2A \log(|\Lambda|)$. Hence,
using the finite speed of propagation result together with Corollary \ref{cor:subex}, we get (details are left to the reader) that
$$
\var_\mu (\PT{t}f) \leq  \frac{||| f |||}{c} e^{-ct} 
+ \frac{1}{c} e^{-c t^{p/(2-p)}} \| f \|_\infty^2
$$
for any $t \geq 1$ and 
for some constant $c$ depending only on 
$p$, $T$, $r$, $d$, $\|W\|_\infty$, $\|W'\|_\infty$, $\|W''\|_\infty$ and $\ell$.
The expected result follows for $t\geq 1$.
Since trivially $\var_\mu(\PT{t}f) \leq \| f - \mu(f) \|_\infty^2$, the expected result follows for any $t>0$.
\end{proof}

%%%%%%%%%%%%%%%%%%%%%%%%%%%%%%%%%%%%%%%%%%%%%%%%%%%%%%%%%%%%%%%%%%%%%%%%%%%%%%%%%%%%%%%%%%%%%%%%%%%%%%%%%
%%%%%%%%%%%%%%%%%%%%%%%%%%%%%%%%%%%%%%%%%%%%%%%%%%%%%%%%%%%%%%%%%%%%%%%%%%%%%%%%%%%%%%%%%%%%%%%%%%%%%%%%%
%%%%%%%%%%%%%%%%%%%%%%%%%%%%%%%%%%%%%%%%%%%%%%%%%%%%%%%%%%%%%%%%%%%%%%%%%%%%%%%%%%%%%%%%%%%%%%%%%%%%%%%%%

\section{The perturbation property improved.} \label{sec:improved}

In this section we prove Proposition \ref{prop:perturbation2} that improves  for cubes the result of Proposition \ref{prop:perturbation}. 
%We 
%a weak Poincar\'e inequality with an almost optimal rate function .
We need to introduce a family of rectangles that will be usefull for our purposes.

%(compare this result with Proposition \ref{prop:perturbation})

Fix $\varepsilon \in(0,1)$.
Let $l_k := (2-\varepsilon)^{k/d}$, and let $\mathbb{F}_k$ be the set of all rectangles $V \Subset \mathbb{Z}^d$
which, modulo translations and permutations of the coordinates, are
contained in
$$
[ 0,l_{k+1} ] \times \dots\times [0,l_{k+d}]
$$

The main property of $\mathbb{F}_k$ is that each rectangle in
$\mathbb{F}_k\setminus \mathbb{F}_{k-1}$ can be obtained as a ``slightly
overlapping union'' of two rectangles in $\mathbb{F}_{k-1}$. More precisely
we have:

\begin{lemma}[\cite{bcc}] \label{lem:geom}
For all $k\in \mathbb{Z}_+$,
for all $\Lambda \in \mathbb{F}_k\setminus \mathbb{F}_{k-1}$ there exists a finite sequence
$\{\Lambda_1^{(i)}, \Lambda_2^{(i)}\}_{i=1}^{s_k}$ in $\mathbb{F}_{k-1}$, where
$s_k := \lfloor l_k^{1/3} \rfloor$, such that, letting $\delta_k := \frac{\varepsilon}{4} \sqrt{l_k}$,
\begin{enumerate}[(i)]
\item $\Lambda = \Lambda_1^{(i)} \cup \Lambda_2^{(i)}$,
\item $d(\Lambda \setminus \Lambda_1^{(i)}, \Lambda \setminus \Lambda_2^{(i)}) \geq \delta_k $,
\item $\left(\Lambda_1^{(i)}\cap \Lambda_2^{(i)}\right) \cap \left(\Lambda_1^{(j)}\cap \Lambda_2^{(j)}\right)
        = \emptyset$,  if $i \neq j$
\end{enumerate}
\end{lemma}

\begin{proof}
The proof is given in \cite[Proposition 3.2]{bcc} for $\varepsilon=1/2$.
The general case given here follows exactly the same line (details are left to the reader).
\end{proof}

\begin{proof}[Proof of Proposition \ref{prop:perturbation2}]
The proof of Proposition \ref{prop:perturbation2} relies on the bisection technique together with
the quasi factorisation of the variance.

The bisection method establishes a simple recursive inequality
between the quantity $\gamma_k(s) := \sup_{\Lambda \in \mathbb{F}_k} \beta_\Lambda(s)$ (recall \eqref{eq:defb}) on
scale $k$ and the same quantity on scale $k-1$. Note that, by construction $\gamma_k$ is non-increasing.

Fix $\Lambda \in \mathbb{F}_k \setminus \mathbb{F}_{k-1}$ and write it as $\Lambda = \Lambda_1 \cup \Lambda_2$ 
with $\Lambda_1,
\Lambda_2\in \mathbb{F}_{k-1}$ satisfying the properties described in Lemma \ref{lem:geom}
above. Without loss of generality we can assume that
all the faces of $\Lambda_1$ and of $\Lambda_2$ lay on the faces of
$\Lambda$ except for one face orthogonal to the first direction $\vec e_1:=(1,0,\cdots,0)$
and that, along that direction, $\Lambda_1$ comes before $\Lambda_2$,
see Figure \ref{fig:bisection}.

\begin{figure}[h]
\psfrag{l}{$\Lambda$}
\psfrag{l1}{$\Lambda_1$}
\psfrag{l2}{$\Lambda_2$}
\psfrag{b}{$\partial_l^r \Lambda_2$}
\includegraphics[width=.50\columnwidth]{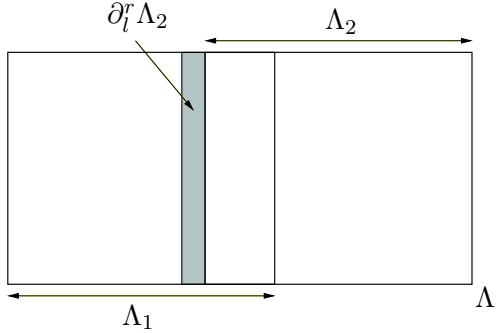}
\caption{The set  $\Lambda = \Lambda_1 \cup \Lambda_2$. The grey region is $\partial_l^r \Lambda_2$.}
\label{fig:bisection}
\end{figure}

We claim that
\begin{claim} There exists $k_0$ such that for $k \geq k_0$,
\begin{equation} \label{eq:claim}
\var_{\mu_\Lambda^\tau} (f) \leq \left(1 + c_1e^{-c_2 \delta_k} \right) \mu_\Lambda^\tau \left( \var_{\mu_{\Lambda_1}}(f) + 
\var_{\mu_{\Lambda_2}}(f) \right)
\end{equation}
for some constant $c_1$ and $c_2$ depending on  $r$ ,$T$, $d$ and $\|W\|_\infty$.
\end{claim}
This bound measures the weak dependence between $\mu_{\Lambda_1}$
and $\mu_{\Lambda_2}$ since it would hold with $c_1=0$ if $\mu_\Lambda^\tau$ was the product 
$\mu_{\Lambda_1} \otimes \mu_{\Lambda_2}$. In other words it is a kind of weak factorisation of the variance.

\begin{proof}[Proof of the Claim]
To prove the claim, let $g$ be a measurable function with respect to $\mathcal{B}_{\Lambda_1^c \cap \Lambda}$. Then, 
by the DLR condition, we have
\begin{equation*}
 \| \mu_{\Lambda_2}(g) - \mu_\Lambda^\tau (g) \|_\infty 
 = 
\| \mu_{\Lambda_2}(g) - \mu_\Lambda^\tau \left( \mu_{\Lambda_2} (g) \right) \|_\infty 
 \leq \!\!
\sup_{\genfrac{}{}{0pt}{}{\eta,\omega \in \Omega :}{\eta_{\Lambda^c}=\omega_{\Lambda^c}} } 
\!\! |\mu_{\Lambda_2}^\eta (g) - \mu_{\Lambda_2}^\omega(g)|
%& \leq &
%C \|g \|_\infty e^{-c d(\Lambda \setminus \Lambda_1^{(i)}, \Lambda \setminus \Lambda_2^{(i)})} .
\end{equation*}
Let $\partial_l^r \Lambda_2 := \{x \in \Lambda \setminus \Lambda_2 \mbox{ such that } x+ i\vec e_1 \in \Lambda_2 \mbox{ for some } i =1,\cdots,r\}$
be the left boundary of width $r$ of $\Lambda_2$, see Figure \ref{fig:bisection}. 
Note that $\sigma \mapsto \mu_{\Lambda_2}^\sigma (g)$ does not depend on
any site $x$ such that $d(x,\Lambda_2)>r$. Hence, if $\eta_{\Lambda^c}=\omega_{\Lambda^c}$, 
$\mu_{\Lambda_2}^\eta (g) - \mu_{\Lambda_2}^\omega(g)$ depends only on the sites in $\partial^r_l \Lambda_2$.
In turn, using a telescopic sum over all $x \in \partial_l^r \Lambda_2$, one has
for any $\eta, \omega \in \Omega$ such that $\eta_{\Lambda^c}=\omega_{\Lambda^c}$,
$$
|\mu_{\Lambda_2}^\eta (g) - \mu_{\Lambda_2}^\omega(g)|
\leq   |\partial_l^r \Lambda_2| 
\sup_{\genfrac{}{}{0pt}{}{x\in \partial_l^r \Lambda_2, \tau,\tau' \in \Omega :}{\tau_{\Lambda_2^c \setminus \{x\}}=\tau'_{\Lambda_2^c \setminus\{x\}}} } 
|\mu_{\Lambda_2}^\tau (g) - \mu_{\Lambda_2}^{\tau'}(g)|.
$$
Now set $h_x:= \frac{Z_{\Lambda_2}^{\tau}}{Z_{\Lambda_2}^{\tau'}} e^{H_{\Lambda_2}^{\tau} - H_{\Lambda_2}^{\tau'}}$
and observe that $h_x$ is a local function with support $\Delta_{h_x}=\{x\}$ and that $\|h_x\|_\infty \leq C$ for some
constant $C=C(r,T,\|W\|_\infty)$.
Then, by a simple  computation and using 
Hypothesis $(H2)$, we have 
$$
|\mu_{\Lambda_2}^\tau (g) - \mu_{\Lambda_2}^{\tau'}(g) |
= 
|\mu_{\Lambda_2}^{\tau} \left(g,h_x\right)|
\leq
C' |\Lambda_1^c \cap \Lambda| \| g \|_\infty e^{- m d (\Lambda \setminus \Lambda_1, \Lambda \setminus \Lambda_2)} 
$$
for some constants $C'$ and $m$ (depending on $r$, $d$, $T$ and $\|W\|_\infty$).
All the previous computations together (recall the definition of $\delta_k$ in Lemma \ref{lem:geom}) lead to
$$
\| \mu_{\Lambda_2}(g) - \mu_\Lambda^\tau (g) \|_\infty \leq C' \| g \|_\infty r\ell_{k+d}^{2d-1} e^{-m \delta_k}
\leq C'' \| g \|_\infty e^{-c_2 \delta_k}
$$
for some constants $C''$ and $c_2$ depending on  $r$ ,$T$, $d$ and $\|W\|_\infty$.

The same holds for $\| \mu_{\Lambda_1}(g) - \mu_\Lambda^\tau (g) \|_\infty$ with
$g$ measurable with respect to $\mathcal{B}_{\Lambda_2^c \cap \Lambda}$. The claim 
follows at once by the following quasi factorisation lemma of \cite{bcc}.
\end{proof}

\begin{lemma}[Quasi factorisation of the Variance \cite{bcc}]
 Let $\Lambda, A, B \Subset \mathbb{Z}^d$ such that $\Lambda=A \cup B$. Assume that for some $\tau \in \Omega$
and $\varepsilon \in [0,\sqrt 2 -1]$,
\begin{eqnarray*}
& \| \mu_B (g) - \mu_\Lambda^\tau(g) \|_\infty \leq  \varepsilon \|g\|_\infty & 
\forall g \in \mathbb{L}^\infty(\Omega,\mathcal{B}_{A^c\cap \Lambda},\mu_\Lambda^\tau) \\
& \| \mu_A (g) - \mu_\Lambda^\tau(g) \|_\infty \leq  \varepsilon \|g\|_\infty & 
\forall g \in \mathbb{L}^\infty(\Omega,\mathcal{B}_{B^c\cap \Lambda},\mu_\Lambda^\tau) .
\end{eqnarray*}
Then,
$$
\var_{\mu_\Lambda^\tau}(f) \leq \frac{1}{1-2\varepsilon - \varepsilon^2}
\mu_\Lambda^\tau \left( \var_{\mu_A}(f) + \var_{\mu_B}(f) \right) \qquad \forall f\in\mathbb{L}^2(\mu_\Lambda^\tau) .
$$
\end{lemma}

\begin{proof}
See \cite[Lemma 3.1]{bcc}
\end{proof}

\begin{remark}
A similar result for the entropy can be found in \cite{cesi}.
\end{remark}

Back to \eqref{eq:claim}, we can use the definition of $\gamma_{k-1}$ twice to get
the following weak Poincar\'e inequality: for any $f : \Omega_\Lambda \to \mathbb{R}$, any $s>0$, it holds
\begin{eqnarray*}
\var_{\mu_\Lambda^\tau}(f)
& \leq & 
\left(1+c_1e^{-c_2 \delta_k} \right) \gamma_{k-1}(s) \left[ \mathcal{D}_\Lambda^\tau (f) + \sum_{x \in \Lambda_1 \cap \Lambda_2} \mu_\Lambda^\tau \left( |\nabla_x f |^2 \right) \right] \\
&& 
\quad + 2 s \left(1+c_1e^{-c_2 \delta_k} \right) \mathrm{Osc}(f) .
\end{eqnarray*}
In order to get rid of the overlapping term $\sum_{x \in \Lambda_1 \cap \Lambda_2} \mu_\Lambda^\tau \left( |\nabla_x f |^2 \right)$ in the latter, as observed in \cite{martinelli}, one can average over the various positions of the pair 
$(\Lambda_1^{(i)},\Lambda_2^{(i)})$ given in Lemma \ref{lem:geom}. In fact, by averaging the previous bound
over the $s_k$ possible choices of $(\Lambda_1^{(i)},\Lambda_2^{(i)})$, we get
\begin{align*}
\var_{\mu_\Lambda^\tau}(f)
& \leq  
\left(1+c_1e^{-c_2 \delta_k} \right) \gamma_{k-1}(s) \left[ \mathcal{D}_\Lambda^\tau (f) + 
\frac{1}{s_k}\sum_{i=1}^{s_k}
\sum_{x \in \Lambda_1^{(i)} \cap \Lambda_2^{(i)}} \!\!\!\! \mu_\Lambda^\tau \left( |\nabla_x f |^2 \right) \right] \\
&
\quad + 2 s \left(1+c_1e^{-c_2 \delta_k} \right) \mathrm{Osc}(f) \\
& \leq 
\left(1+c_1e^{-c_2 \delta_k} \right) \left(1+\frac{2}{s_k}\right) \gamma_{k-1}(s) \mathcal{D}_\Lambda^\tau (f) \\
&
\quad 
+ 2 s \left(1+c_1e^{-c_2 \delta_k} \right) \mathrm{Osc}(f) .
\end{align*}
In the last line we used that $\left(\Lambda_1^{(i)}\cap \Lambda_2^{(i)}\right) \cap \left(\Lambda_1^{(j)}\cap \Lambda_2^{(j)}\right) = \emptyset$  for $i \neq j$, {\it i.e.} Point $(iii)$ of Lemma \ref{lem:geom}.
It follows that
$$
\gamma_k(s) \leq 
\left(1+c_1e^{-c_2 \delta_k} \right) \left(1+\frac{2}{s_k}\right) \gamma_{k-1} \left( \frac{s}{2 \left(1+c_1e^{-c_2 \delta_k} \right)} \right) \qquad \forall s>0 .
$$
By iteration, we get for any $k \geq k_0$ and any $s>0$,
$$
\gamma_k(s) \leq 
\prod_{i=k_0+1}^k \!\!  \left(1+c_1e^{-c_2 \delta_i} \right) \left(1+\frac{2}{s_i}\right) 
\gamma_{k_0} \left( \frac{s}{2^{k-k_0} \prod_{i=k_0+1}^k \left(1+c_1e^{-c_2 \delta_i} \right)} \right) \! .
$$
Note that for some $C=C(r,T,d,\|W\|_\infty)$,
$$
1 \leq \prod_{i=k_0+1}^k \left(1+c_1e^{-c_2 \delta_i} \right) \leq 
\prod_{i=0}^\infty \left(1+c_1e^{-c_2 \delta_i} \right) \leq C
$$
and similarly for  $\prod_{i=k_0+1}^k \left(1+\frac{2}{s_i}\right)$. Hence, since $\gamma_{k_0}$ is non-increasing,
$$
\gamma_k(s) \leq C^2 \gamma_{k_0}\left( \frac{s}{C2^k} \right) \qquad \forall k \geq k_0, \;\forall s>0 .
$$
We are left with an estimate of $\gamma_{k_0}$. This is given by Proposition \ref{prop:perturbation}. Indeed,
since $k_0$ is a constant depending only on $r$ ,$T$, $d$ and $\|W\|_\infty$, 
Proposition \ref{prop:perturbation} guarantees that $\gamma_{k_0}(s) \leq C' \beta \left(s/C' \right)$
for some $C'=C'(r,T,d,\|W\|_\infty)$.

In conclusion, we have proved that for any $\Lambda \in \mathbb{F}_k$,
$$
\var_{\mu_\Lambda^\tau} (f) \leq C'' \beta \left( \frac{s}{2^k C''} \right) \mathcal{D}_\Lambda^\tau (f) 
+ s \mathrm{Osc}(f) \qquad \forall \tau \in \Omega, \; \forall f, \; \forall s>0
$$
for some $C''=C''(r,T,d,\|W\|_\infty)$. 

Now consider a volume $\Lambda=[-L,L]^d$. Observe that $\Lambda \in \mathbb{F}_k$ as soon as
$2L \leq l_{k+1}$. Take $k$ to be the smallest integer satisfying such a property.
After some computations, this leads to $2^k \leq c |\Lambda|^{\frac{\log 2}{\log(2-\varepsilon)}}$
for some universal constant $c>0$. Since $\beta_\Lambda$ is non-increasing, we get the expected result.
This achieves the proof.
\end{proof}

%%%%%%%%%%%%%%%%%%%%%%%%%%%%%%%%%%%%%%%%%%%%%%%%%%%%%%%%%%%%%%%%%%%%%%%%%%%%%%%%%%%%%%%%%%%%%%%%%%%%%%%%%
%%%%%%%%%%%%%%%%%%%%%%%%%%%%%%%%%%%%%%%%%%%%%%%%%%%%%%%%%%%%%%%%%%%%%%%%%%%%%%%%%%%%%%%%%%%%%%%%%%%%%%%%%
%%%%%%%%%%%%%%%%%%%%%%%%%%%%%%%%%%%%%%%%%%%%%%%%%%%%%%%%%%%%%%%%%%%%%%%%%%%%%%%%%%%%%%%%%%%%%%%%%%%%%%%%%

\section{Finite speed of propagation} \label{sec:fs}

This section is dedicated to the proof of Proposition \ref{prop:finitespeed} (we recall below) on the finite speed of propagation. This result is somehow standard and would certainly not surprise the specialits. 
Nevertheless we give the proof for completeness.

Recall the definition of the finite volume and infinite volume 
Markov semi-groups $(\PT{t}^{\Lambda,\tau})_{t \geq 0}$ and $(\PT{t})_{t\geq 0}$. Recall also the definition of $|||f|||$.

\begin{proposition}[Finite speed of propagation] 
Assume $(H_1)$ and $(H2)$. Fix and integer $\ell \geq 1$.
Then, for any local function $f$ with support $\Delta_f \subset [-\ell,\ell]^d$, 
any $L$ multiple of $r$, any boundary condition $\tau \in \Omega$,
$$
\| \PT{t}f - \PT{t}^{\Lambda,\tau} f \|_\infty \leq C ||| f ||| \left( \frac{C't}{L} \right)^{C''L} e^{Ct} 
\qquad \forall t>0
$$
with $\Lambda = [-L,L]^d$,
for some constant $C, C', C''>0$ depending only on $r$, $d$, $\|W'\|_\infty$, $\|W''\|_\infty$ and $\ell$.
\end{proposition}

\begin{remark}
Note that this bound is particularly interesting when $L \gg t$.
\end{remark}

\begin{proof}
We follow \cite{zegarlinski}.
Fix $t>0$, $\Lambda=[-L,L]^d \Subset \mathbb{Z}^d$ with $L$ a multiple of $r$,
a boundary condition $\tau \in \Omega$
and a local function $f$ with support $\Delta_f$ 
containing $0$. Then,
\begin{equation}\label{eq:10}
(\PT{t}-\PT{t}^{\Lambda,\tau})f = - \int_0^t \left(\frac{d}{ds} (\PT{t-s} \PT{s}^{\Lambda,\tau}) f\right) ds
= \int_0^t \PT{t-s}(\GI - \GI_\Lambda^\tau)\PT{s}^{\Lambda,\tau} f ds .
\end{equation}
For simplicity let $f_s^\Lambda:= \PT{s}^{\Lambda,\tau}f$ and note that 
its support $\Delta_{f_s^\Lambda} \subset \Lambda$. 
Therefore, 
\begin{eqnarray}
(\GI - \GI_\Lambda^\tau) f_s^\Lambda 
& =  &
\sum_{x \in \Lambda} \left( \nabla_x H_\Lambda - \nabla_x H_\Lambda^\tau \right) \cdot \nabla_x f_s^\Lambda \nonumber  \\
& = &
 \sum_{x \in \Lambda: d(x,\Lambda^c)\leq r} 
\left( \nabla_x H_\Lambda - \nabla_x H_\Lambda^\tau \right) \cdot \nabla_x f_s^\Lambda .
\end{eqnarray}
Now our aim is to control $\nabla_x f_s^\Lambda$.

Take $y \in \Lambda$. By definition of $L_\Lambda^\tau$, we have
$$
[\nabla_y,L_\Lambda^\tau] := 
\nabla_y L_\Lambda^\tau - L_\Lambda^\tau \nabla_y = 
\sum_{x \in \Lambda} \nabla_y \nabla_x H_\Lambda^\tau \cdot \nabla_x
=
\!\!\!\! \sum_{x \in \Lambda:d(x,y)\leq r} \!\!\!\! \nabla_y \nabla_x H_\Lambda^\tau \cdot \nabla_x .
$$
Thus, (we skip the superscrip $\tau$)
\begin{eqnarray*}
\nabla_y f_s^\Lambda & = & 
\PT{s}^{\Lambda} \nabla_y f + 
\int_0^s \left( \frac{d}{du} \PT{s-u}^{\Lambda} \nabla_y f_u^\Lambda \right) du \\
& = & 
\PT{s}^{\Lambda} \nabla_y f +  \int_0^s \PT{s-u}^{\Lambda} [\nabla_y,L_\Lambda^\tau] f_u^\Lambda du \\
& = & 
\PT{s}^{\Lambda} \nabla_y f + \!\!\!\! \sum_{x \in \Lambda:d(x,y)\leq r}  \int_0^s \PT{s-u}^{\Lambda}
\nabla_y \nabla_x H_\Lambda^\tau \cdot \nabla_x f_u^\Lambda du .
\end{eqnarray*}
Hence, thanks to Hypothesis $(H_1)$ and the fact that $\PT{t}^\Lambda$ is a contraction in the sup norm,
\begin{equation}\label{eq:step1}
\|\nabla_y f_s^\Lambda \|_\infty  
\leq 
\| \nabla_y f \|_\infty + \|W''\|_\infty \sum_{x \in \Lambda:d(x,y)\leq r}
\int_0^s  \| \nabla_x f_u^\Lambda \|_\infty du \;.
\end{equation}
Then, for any $ n =0,1,\ldots, L/r$, define 
$$
Y_n(u) := \sum_{x \in \Lambda : d(x,\Lambda^c) \leq L-nr} \| \nabla_x f_u^\Lambda \|_\infty .
$$
Recall that $\Delta_f \subset [-\ell,\ell]^d$.
Since $\nabla_x f = 0$ unless $x \in \Delta_f$, we get from
\eqref{eq:step1} that for $n=\ell+1,\ldots, L/r$,
$$
Y_n(s) \leq (2r)^d \|W''\|_\infty \int_0^s Y_{n-1}(u) du \;.
$$
On the other hand, for $n=0,1,\ldots, \ell$,
$$
Y_n(s) \leq |||f||| + (2r)^d \|W''\|_\infty \int_0^s Y_{n-1}(u) du,
$$
with the convention that $Y_{-1}:=Y_0$.

It follows that $Y_n(t) \leq |||f||| \exp \{ C t \}$ for any
$0 \leq n \leq \ell$, with 
$C:=(2r)^d \|W''\|_\infty$. Moreover, an easy induction gives for any $\ell < n \leq L/r$
\begin{equation} \label{eq:b}
Y_n(t) \leq  |||f||| R(n-\ell,t) 
\qquad \mbox{with }
R(m,t):= e^{Ct} - \sum_{k=0}^{m} \frac{(Ct)^k}{k!} 
\leq \left(\frac{Cte}{m} \right)^m  e^{Ct}.
\end{equation}
Finally, using the fact that $\PT{t}$ is a contraction in the sup norm and Hypothesis $(H1)$,
we get from \eqref{eq:10} and \eqref{eq:b} that
\begin{eqnarray*}
\| (\PT{t}-\PT{t}^{\Lambda,\tau})f \|_\infty
& \leq  &
 \int_0^t \| (\GI - \GI_\Lambda^\tau)f_s^\Lambda \|_\infty ds  \\
& \leq & 
\sum_{x \in \Lambda : d(x,\Lambda^c) \leq r} 
\int_0^t \| \nabla_x H_\Lambda - \nabla_x H_\Lambda^\tau \|_\infty  \|f_s^\Lambda \|_\infty ds \\
& \leq & 
4(2r)^d \|W'\|_\infty \int_0^t Y_{\frac{L}{r}-1} (s) ds \\
& \leq &
\frac{4(2r)^d \|W'\|_\infty}{C} |||f||| R\left(\frac{L}{r}-\ell,t\right) \\
\end{eqnarray*}
The expected result follows from \eqref{eq:b}. This achieves the proof.
\end{proof}

\bibliographystyle{plain}
\bibliography{bib-weak-poincare-mecastat}

\end{document}